\documentclass[11pt]{article}
\usepackage[T2A]{fontenc}
\usepackage[cp1251]{inputenc}
\usepackage{amsmath}
\usepackage{amssymb}
\def\a{\mathbf{a}}

\def\m{\mathbf{m}}
\def\f{\mathbf{f}}

\def\b{\mathbf{b}}
\def\p{\mathbf{p}}
\def\N{\mathbb{N}}
\def\Z{\mathbb{Z}}

\def\O{\mathcal{O}}

\def\zetta{\boldsymbol{\zeta}}
\def\betta{\boldsymbol{\beta}}

\def\aalpha{\boldsymbol{\alpha}}
\def\betta{\boldsymbol{\beta}}

\newtheorem{theorem}{\hspace*{\parindent}Theorem}

\title{Extensions of Karlsson-Minton summation theorem and some consequences of the first Miller-Paris transformation}
\author{D.B.\:Karp$^{\rm a,b}$\footnote{Corresponding author. E-mail: D. Karp -- \emph{dimkrp@gmail.com},
E.\:Prilepkina --  \emph{pril-elena@yandex.ru}}~~and
E.G.\:Prilepkina$^{\rm a,b}$
\\[10pt]
\\
\small{\textit{$\phantom{1}^a$Far Eastern Federal University, 8
Sukhanova street, Vladivostok, 690950, Russia}}
\\
\small{\textit{$\phantom{1}^b$Institute of Applied Mathematics,
FEBRAS, 7 Radio Street, Vladivostok,  690041, Russia}}}
\date{}

\begin{document}
\maketitle

\begin{center}
\parbox{12cm}{
\small\textbf{Abstract.}
In this paper we give several independent extensions of the Karlsson-Minton summation formula for the generalized hypergeometric function with integral parameters differences. In particular, we examine the ''prohibited'' values for the integer top parameter in Minton's formula, extend one unit negative difference in Karlsson's formula to a finite number of integer negative  differences and establish known and new summation and transformation formulas when the unit negative difference is allowed to take arbitrary values.  We also present a recurrence relation reducing the case of integer negative difference to the Karlsson-Minton case of unit negative difference. Further, we explore some alternative forms of the first Miller-Paris transformation, including one expressed in terms of Meijer-N{\o}rlund $G$ function.}
\end{center}

\bigskip


\bigskip


\bigskip

\section{Introduction}
Throughout the paper we will use the standard  symbols $\Z$ and $\N$ to denote the integer and natural numbers, respectively; $\N_0=\N\cup\{0\}$.
Suppose  $\a=(a_1,\ldots,a_p)$ and $\b=(b_1,\ldots,b_q)$ are complex vectors, $(a)_0=1$, $(a)_n=a(a+1)\ldots(a+n-1)$, $n\in\N$,  is the Pochhammer's symbol. We will use the following convenient abbreviated notation:
$$
\Gamma(\a)=\Gamma(a_1)\Gamma(a_2)\cdots\Gamma(a_p),~~(\a)_n=(a_1)_n(a_2)_n\cdots(a_p)_n,~~\a+\mu=(a_1+\mu,a_2+\mu,\dots,a_p+\mu);
$$
$\a_{[k]}$ will denote the vector $\a$ with $k$-th component removed. Inequalities like  $\Re(\a)>0$ and properties like  $-\a\notin\N$ will be understood element-wise, i.e. $\Re(a_i)>0$, $-a_i\notin\N$ for $i=1,\ldots,p$ in the above example. The generalized hypergeometric functions is defined by the series
\begin{equation}\label{eq:hypergeometic}
{}_{p}F_{q}\left.\!\!\left(\begin{matrix}\a\\\b\end{matrix}\right\vert x\right)=\sum\limits_{n=0}^{\infty}\frac{(\a)_n x^n}{(\b)_nn!}
\end{equation}
for the values of $x$ making the series convergent and by analytic continuation elsewhere.  The size of a vector appearing in the parameters of the generalized hypergeometric function can usually be read off its indices.  More detailed information about the definition and properties of the generalized hypergeometric function can be found in standard textbooks \cite[Section~2.1]{AAR}, \cite[Chapter~12]{BealsWong}) and reference books \cite[Section~5.1]{LukeBook}, \cite[Sections 16.2-16.12]{NIST}. The argument $x=1$  of the generalized hypergeometric function will be routinely omitted: ${_{p}F_q}(\a;\b):={_{p}F_q}(\a;\b;1)$.
In what follows we also will use the notation $\m=(m_1,\ldots,m_r)\in\N^r$, $m=m_1+m_2+\ldots+m_r$, $\f=(f_1,\ldots,f_r)$ and $(\f)_{\m}=(f_1)_{m_1}\cdots(f_r)_{m_r}$.

In 1970 Minton \cite{Minton} proved the summation formula
\begin{equation}\label{eq:Minton}
{}_{r+2}F_{r+1}\left.\!\!\left(\begin{matrix}-k,b,\f+\m\\b+1,\f\end{matrix}\right.\right)=\frac{k!}{(b+1)_{k}}
\frac{(\f-b)_\m}{(\f)_\m}
\end{equation}
valid for $k\ge{m}$, $k\in\N$.  Soon thereafter, his  result was generalized by Karlsson \cite{Karlson} who replaced $-k$ by an arbitrary complex number $a$ satisfying $\Re(1-a-m)>0$ to get
 \begin{equation}\label{eq:Karl1}
{}_{r+2}F_{r+1}\left.\!\!\left(\begin{matrix}a,b,\f+\m\\b+1,\f\end{matrix}\right.\right)=\frac{\Gamma(b+1)\Gamma(1-a)}
{\Gamma(b+1-a)}\frac{(\f-b)_\m}{(\f)_\m}.
\end{equation}
Gasper \cite{Gasper} deduced a $q$-analogue and a generalization of Minton's and Karlsson's formulas; Chu \cite{Chu1,Chu2} found extensions to bilateral hypergeometric and $q$-hypergeometric series; their results were re-derived by simpler means and further generalized by Schlosser \cite{Schlosser2}, who also found multidimensional extensions  to hypergeometric functions associated with root systems \cite{Schlosser1}.  For further developments in this directions, see also \cite{Rosengren1,Rosengren2}.  Quite recently, Miller and Srivastava \cite{MS2010} found some alternative derivations of the original Karlsson and Minton formulas.  Some of the results mentioned above found applications in the theory of rook polynomials and related combinatorial problems as discovered by Haglund in \cite{Haglund}, who also noticed that some of the results can be obtained by specializing a very general transformation due to Slater.  We also mention an interesting work \cite{LVW} by Letessier, Valent  and Wimp, where the reduction of the differential equation satisfied by the generalized hypergeometric functions with some integer parameter differences was found.

Conditions  $k\ge{m}$ in (\ref{eq:Minton}) and $\Re(1-a-m)>0$ in (\ref{eq:Karl1}) can not be dispensed with: the formulas
become false when they are violated (take for example $a=-4$, $b=33/17$, $c=b+1$, $f_1=4.2$, $f_2=-5/3$, $m_1=7$, $m_2=8$).
The situation with (\ref{eq:Karl1}) is thereby quite subtle:  the function on the right hand side is continuous at $a=1-m$ while the function on the left hand side is finite but has a discontinuity at this point. Further, if $\Re(1-a-m)<0$ and $a$ is not a negative integer or zero, the series on the left hand side diverges. At the same time, it reduces to a finite sum for $a=-k$, $k\in\N$, so that such points are ''isolated summability points'' and analytic continuation of the right hand side in $a$ does not lead to the correct value of the left hand side at such point. In section~2 of this paper we find the correct value of the left hand side of Minton's formula (\ref{eq:Minton}) at these special points $k=0,1,\ldots,m-1$.

In the same section we present an extension of Karlsson's formula (\ref{eq:Karl1}) obtained by replacing the number $b$ by a vector $\b$ and the number $b+1$ by the vector $\b+\p$, where $\p$ is a positive integer vector. Note that the Clausen series ${}_3F_2$ with two negative integer parameter differences has been studied recently by Shpot and Srivastava in \cite{ShpotSriv}. Next we present an alternative proof of another generalization of Karlsson's formula (\ref{eq:Karl1}) due to Gasper \cite[(18)]{Gasper}. This generalization has the form of a transformation formula for the generalized hypergeometric function on the left hand side of (\ref{eq:Karl1}) with $b+1$ replaced with an arbitrary number $c$ (restricted to ensure convergence), of which Karlsson's formula is a particular case. For our proof we will need the generalized Stieltjes representation of the generalized hypergeometric function.  More precisely, if $\min_i\{\Re a_i\}>0$ and $\Re\{\sum_{i=1}^q(b_i-a_i)\}>0$, then \cite[(2)]{KLJAT2017}
\begin{equation}\label{eq:int2}
{}_{q+1}F_{q}\left.\!\!\left(\begin{matrix}\sigma,\a\\\b\end{matrix}\right.\right)=
\frac{\Gamma(\b)}{\Gamma(\a)}\int\limits_0^{1}G_{q,q}^{q,0}\!\left(\!t~\vline\begin{array}{l}\b-1\\
\a-1\end{array}\!\!\right)\frac{dt}{(1-t)^\sigma}.
\end{equation}
Here $G_{q,q}^{q,0}$ is the Meijer-N{\o}rlund $G$ function defined by the Mellin-Barnes integral of the form
\begin{equation}\label{eq:G-defined}
G^{q,0}_{q,q}\!\left(\!z~\vline\begin{array}{l}\b\\\a\end{array}\!\!\right)\!\!:=
\\
\frac{1}{2\pi{i}}
\int\limits_{\mathcal{L}}\!\!\frac{\Gamma(\a\!+\!s)}{\Gamma(\b+\!s)}z^{-s}ds.
\end{equation}
We omit the details regarding the choice of the contour ${\mathcal{L}}$ as the definition of (the general case of) Meijer's $G$ function can be found in standard text- and reference- books \cite[section~5.2]{LukeBook}, \cite[16.17]{NIST}, \cite[8.2]{PBM3} and \cite[Chapter~12]{BealsWong}.  See, also our papers \cite{KLJAT2017,KPSIGMA}.  The next shifting property is straightforward from the definition (\ref{eq:G-defined}), but nevertheless is very useful \cite[8.2.2.15]{PBM3}:
\begin{equation}\label{eq:Gtimespower}
z^{\alpha}G^{q,0}_{q,q}\!\left(\!z~\vline\begin{array}{l}\b\\\a\end{array}\!\!\right)=G^{q,0}_{q,q}\!\left(\!z~\vline\begin{array}{l}\b+\alpha\\\a+\alpha\end{array}\!\!\right),
\end{equation}
where $\alpha$  is any complex number.  We further mention that the coefficients $g_n(\a;\b)$ of the power series expansion
\begin{equation}\label{eq:Norlund}
G^{q,0}_{q,q}\!\left(\!1-z~\vline\begin{array}{l}\b\\\a,0\end{array}\!\!\right)=\frac{z^{\nu_q-1}}{\Gamma(\nu_q)}
\sum\limits_{n=0}^{\infty}\frac{g_n(\a;\b)}{(\nu_q)_n}z^n,
\end{equation}
where $\nu_q=\nu_q(\a;\b)=\sum_{j=1}^{q}b_j-\sum_{j=1}^{q-1}a_j$, are known as N{\o}rlund's coefficients (see \cite[(11)]{KLJAT2017}).

Among the most important discoveries in the field of hypergeometric transformations over the last decade are
transformation formulas for the generalized hypergeometric functions found, mainly, by Miller and Paris (with important contributions by Rathie, Kim and others) and summarized in the seminal paper \cite{MP2013}. In section~3 of this paper we explore some consequences of the Miller-Paris transformations, including a guise of the first transformation expressed as a transformation for the Meijer-N{\o}rlund $G$ function,  an identity  for the products of the generalized hypergeometric functions at unit argument and a transformation formula for right hand side of the identity
\begin{equation}\label{eq:sumMil}
{}_{r+2}F_{r+1}\left.\!\!\left(\begin{matrix}a,b,\f+\m\\c,\f\end{matrix}\right.\right)=
\frac{\Gamma(c)\Gamma(c-a-b)}{\Gamma(c-a)\Gamma(c-b)}\sum\limits_{k=0}^m
{}_{r+1}F_{r}\left.\!\!\left(\begin{matrix}-k,\f+\m\\\f\end{matrix}\right.\right)
\frac{(a)_k(b)_k}{(1+a+b-c)_k k!}
\end{equation}
demonstrated in \cite[Theorem~4]{MP2012}. Finally, we find a recurrence relation for the generalized hypergeometric
function with integer shifts.

\section{Extensions of Karlsson-Minton theorem}

To formulate our first theorem  we will need an explicit formula for N{\o}rlund's coefficients defined in (\ref{eq:Norlund}). Assuming, as before, that $\a=(a_1,\ldots,a_{q-1})$ and $\b=(b_1,\ldots,b_q)$ are some complex vectors, we recall that the coefficients $g_n(\a;\b)$ are symmetric polynomials in the components of $\a$ and $\b$ (separately) given explicitly by \cite[(2.9)]{KPSIGMA}
\begin{equation}\label{eq:Norlund-explicit}
g_n(\a;\b)=\sum\limits_{0\leq{j_{1}}\leq{j_{2}}\leq\cdots\leq{j_{q-2}}\leq{n}}
\prod\limits_{l=1}^{q-1}\frac{(\psi_l+j_{l-1})_{j_{l}-j_{l-1}}}{(j_{l}-j_{l-1})!}(b_{l+1}-a_{l})_{j_{l}-j_{l-1}},
\end{equation}
where $\psi_l=\sum_{i=1}^{l}(b_i-a_i)$, $j_0=0$, $j_{q-1}=n$. Their further properties can be found in our recent papers \cite{KLJAT2017,KPSIGMA}. In particular, we showed in \cite[Theorem~2]{KPSIGMA} that the initial coefficients are given by
\begin{equation}\label{eq:Norlund18}
\begin{split}
&g_0(\a;\b)=1,~~~~g_1(\a;\b)=\sum_{l=1}^{q-1}(b_{l+1}-a_l)\psi_l,
\\
&g_2(\a;\b)=\frac{1}{2}\sum_{l=1}^{q-1}(b_{l+1}-a_l)_2(\psi_l)_2+\sum_{k=2}^{q-1}(b_{k+1}-a_{k})(\psi_{k}+1)\sum_{l=1}^{k-1}(b_{l+1}-a_l)\psi_l.
\end{split}
\end{equation}
Below, we extend Minton's formula (\ref{eq:Minton}) to the values of $k\in\{0,1,\ldots,m-1\}$.

\begin{theorem}\label{thm:Karlson5}
Suppose that $k\in\N_0$, $0\leq k\leq m-1$, $(b+1)_k\ne0$, $(\f)_{\m}\ne0$,  Then
\begin{equation}\label{eq:Karlson12}
{}_{r+2}F_{r+1}\left.\!\!\left(\begin{matrix}-k,b,\f+\m\\b+1,\f\end{matrix}\right.\right)=
\frac{k!}{(b+1)_{k}}
\frac{(\f-b)_\m}{(\f)_\m}-\frac{(-1)^mk!b}{(\f)_\m}q_k,
\end{equation}
where $q_k=\sum\limits_{i=0}^{m-k-1}g_{m-k-i-1}(\aalpha;\betta)(b-i)_i$, $\aalpha=(b-\f),$ $\betta=(b-\f-\m,b+k)$ and the coefficients $g_n$ are defined in \emph{(\ref{eq:Norlund-explicit})}.
In particular,
\begin{equation}\label{eq:Karlson14}
{}_{r+2}F_{r+1}\left.\!\!\left(\begin{matrix}1-m,b,\f+\m\\b+1,\f\end{matrix}\right.\right)=\frac{(m-1)!}{(b+1)_{m-1}}
\frac{(\f-b)_\m}{(\f)_\m}
-\frac{(-1)^{m}(m-1)!b}{(\f)_\m}
\end{equation}
and
\begin{multline}\label{eq:Karlson15}
{}_{r+2}F_{r+1}\left.\!\!\left(\begin{matrix}2-m,b,\f+\m\\b+1,\f\end{matrix}\right.\right)
=\frac{(m-2)!}{(b+1)_{m-2}}\frac{(\f-b)_\m}{(\f)_\m}
\\
+\frac{(-1)^{m}(m-2)!b}{(\f)_\m}\left(1-b+m(m-2+f_r)+
\sum\limits_{q=1}^{r-1}(f_q-f_{q+1}-m_{q+1})\sum\limits_{i=1}^q m_i\right).
\end{multline}
\end{theorem}
\textbf{Proof.} According to \cite[(2.4)]{KPSIGMA} when $f_i-f_j\notin\Z$ and $f_i-c\notin\Z$ the next identity holds:
\begin{equation}\label{eq:Ghiper}
G_{r+1,r+1}^{r+1,0}\!\left(\!z~\vline\begin{array}{l}c,\f\\
b,\f+\m\end{array}\!\!\right)=z^b\frac{\Gamma(\f+\m-b)}{\Gamma(c-b)\Gamma(\f-b)}
{}_{r+1}F_{r}\left.\!\!\left(\begin{matrix}1-c+b,1-\f+b\\1-\f-\m+b\end{matrix}\right\vert
z\right).
\end{equation}
In view of the shifting property (\ref{eq:Gtimespower}), this implies that
\begin{equation}\label{eq:Ghiper1}
G_{r+1,r+1}^{r+1,0}\!\left(\!t~\vline\begin{array}{l}b+k,b-\f-\m\\
b-1,b-\f\end{array}\!\!\right)=t^{b-1}\frac{(-1)^{m}(\f)_\m}{k!}
{}_{r+1}F_{r}\left.\!\!\left(\begin{matrix}-k,\f+\m\\\f\end{matrix}\right\vert
t\right).
\end{equation}
Integrating this formula with respect to $t$ from $0$ to $1$ we get on the left hand side:
\begin{multline}\label{eq:intGhiper1}
\int\limits_{0}^1G_{r+1,r+1}^{r+1,0}\!\left(\!t~\vline\begin{array}{l}b+k,b-\f-\m\\
b-1,b-\f\end{array}\!\!\right)dt
=\frac{(-1)^{m}(\f)_\m}{k!}\sum\limits_{n=0}^\infty\frac{(-k)_n(\f+\m)_n}{n!(\f)_n}\int\limits_0^1t^{b-1+n}dt
\\
=\frac{(-1)^{m}(\f)_\m}{k!}\sum\limits_{n=0}^\infty\frac{(-k)_n(\f+\m)_n}{n!(\f)_n(b+n)}
=\frac{(-1)^{m}(\f)_\m}{k!b}
{}_{r+2}F_{r+1}\left.\!\!\left(\begin{matrix}-k,b,\f+\m\\b+1,\f\end{matrix}\right.\right).
\end{multline}
On the other hand, it was shown by N{\o}rlund that \cite[(17)]{KLJAT2017}
\begin{equation}
\int\limits_{0}^1G_{r+1,r+1}^{r+1,0}\!\left(\!t~\vline\begin{array}{l}b+k,b-\f-\m\\
b-1,b-\f\end{array}\!\!\right)dt=\frac{(1-\f-\m+b)_\m}{(b)_{k+1}}-q_k,
\end{equation}
which yields (\ref{eq:Karlson12}) once we account for the straightforward identity
\begin{equation}{\label{eq:Pocham}}
\frac{(\f-b)_\m(1-\f-\m)_\m}{(\f)_\m(1-\f+b-\m)_\m}=1.
\end{equation}
Using explicit formulas (\ref{eq:Norlund18}) for N{\o}rlund's coefficients for $k=m-1$ and $k=m-2$ we obtain (\ref{eq:Karlson14})
and  (\ref{eq:Karlson15}). Additional assumptions made in the course of the proof can now be removed by continuity of both sides of (\ref{eq:Karlson12}).
$\hfill\square$

\begin{theorem}\label{thm:Karlson8}
Suppose $\b=(b_1,\ldots,b_l)$ is a complex vector,  $\p=(p_1,\ldots,p_l)$ is a vector of positive integers, $p=p_1+p_2+\ldots+p_l$, and all elements of the vector $\betta=(b_1,b_1+1,\ldots,b_1+p_1-1,\ldots,b_l,b_l+1,\ldots,b_l+p_l-1)=(\beta_1,\beta_2,\ldots,\beta_p)$ are distinct. If $\Re(p-a-m)>0$, then
\begin{equation}\label{eq:Karlson9}
\frac{1}{\Gamma(1-a)}{}_{r+l+1}F_{r+l}\left.\!\!\left(\begin{matrix}a,\b,\f+\m\\\b+\p,\f\end{matrix}\right.\right)
=\frac{(\b)_\p}{(\f)_\m}\sum\limits_{q=1}^{p}\frac{\Gamma(\beta_q)(\f-\beta_q)_\m}{B_q\Gamma(1+\beta_q-a)},
\end{equation}
where $B_q=\prod\limits_{v=1,\,v\neq q}^{p}(\beta_v-\beta_q)$.

In particular, if $p\in\N$ and $\Re(p-a-m)>0$, then
\begin{equation}\label{eq:Karlson3}
\frac{1}{\Gamma(1-a)}{}_{r+2}F_{r+1}\left.\!\!\left(\begin{matrix}a,b,\f+\m\\b+p,\f\end{matrix}\right.\right)
=\frac{\Gamma(b+p)}{(\f)_\m}\sum\limits_{q=0}^{p-1}\frac{(b)_q(\f-b-q)_\m}{(p-q-1)!(-q)_q\Gamma(1+q+b-a)}.
\end{equation}
\end{theorem}
\textbf{Proof.} Using the  the partial fraction decomposition
$$
\frac{1}{(\betta+x)_1}=\sum\limits_{q=1}^{p}\frac{B_q^{-1}}{\beta_q+x},
$$
with $B_q=\prod\limits_{v=1,\,v\neq q}^{p}(\beta_v-\beta_q)$, we get
$$
\frac{(\b)_n}{(\b+\p)_n}=\frac{(\b)_\p}{(\b+n)_\p}=\frac{(\b)_{\p}}{(\betta+n)_1}=(\b)_\p\sum\limits_{q=1}^{p}\frac{B_q^{-1}(\beta_q)_n}{(\beta_q+1)_n(\beta_q)_1}.
$$
Hence,
\begin{equation}\label{eq:Karlson10}
{}_{r+l+1}F_{r+l}\left.\!\!\left(\begin{matrix}a,\b,\f+\m\\\b+\p,\f\end{matrix}\right.\right)
=\sum\limits_{q=1}^{p}\frac{(\b)_\p}{B_q(\beta_q)_1}{}_{r+2}F_{r+1}\left.\!\!\left(\begin{matrix}a,\beta_q,\f+\m\\\beta_q+1,\f\end{matrix}\right.\right).
\end{equation}
It remains to apply (\ref{eq:Karl1}) to obtain (\ref{eq:Karlson9}) under additional restriction $\Re(1-a-m)>0$. This restriction can now be relaxed to
$\Re(p-a-m)>0$ by appealing to analytic continuation. $\hfill\square$

\medskip

For  $p=2$  after elementary manipulations the summation formula (\ref{eq:Karlson3}) takes the form:
 \begin{equation}\label{eq:Karlson2}
{}_{r+2}F_{r+1}\left.\!\!\left(\begin{matrix}a,b,\f+\m\\b+2,\f\end{matrix}\right.\right)
=-\frac{\Gamma(b+2)\Gamma(1-a)\left(b(\f-1-b)_\m+(a-b-1)(\f-b)_\m\right)}{\Gamma(b+2-a)(\f)_\m}.
\end{equation}

As another example, take $\b=(b,c)$,  $\p=(1,1)$ in Theorem~\ref{thm:Karlson8}  to get:
$$
{}_{r+3}F_{r+2}\left.\!\!\left(\begin{matrix}a,b,c,\f+\m\\b+1,c+1,\f\end{matrix}\right.\right)
=\frac{\Gamma(1-a)bc}{(c-b)(\f)_\m}\left(\frac{\Gamma(b)(\f-b)_\m}{\Gamma(1+b-a)}-\frac{\Gamma(c)(\f-c)_\m}{\Gamma(1+c-a)}\right).
$$

In \cite[(18)]{Gasper} Gasper discovered a transformation formula which extends Karlsson's summation theorem if $b+1$ in (\ref{eq:Karl1}) is replaced by an indeterminate $c$.
In the following theorem we present an alternative derivation of Gasper's formula (\ref{eq:Kar4}) based on the integral representation (\ref{eq:int2}).
\begin{theorem}\label{thm:Karlson4}
If $\Re(c-a-b-m)>0$, $-c\notin\N_0$ and $-f_i\notin\N_0$, then
\begin{equation}\label{eq:Kar4}
\frac{1}{\Gamma(1-a)}{}_{r+2}F_{r+1}\left.\!\!\left(\begin{matrix}a,b,\f+\m\\c,\f\end{matrix}\right.\right)
=
\frac{\Gamma(c)(\f-b)_\m}{\Gamma(1+b-a)\Gamma(c-b)(\f)_\m}
{}_{r+2}F_{r+1}\left.\!\!\left(\begin{matrix}b,1-c+b,1-\f+b\\1+b-a,1-\f-\m+b\end{matrix}\right.\right).
\end{equation}
\end{theorem}
\textbf{Proof.}  Indeed, combining (\ref{eq:int2}) and (\ref{eq:Ghiper}) we get
 \begin{multline}
{}_{r+2}F_{r+1}\left.\!\!\left(\begin{matrix}a,b,\f+\m\\c,\f\end{matrix}\right.\right)=
\frac{\Gamma(c)(\f-b)_\m}{\Gamma(b)\Gamma(c-b)(\f)_\m}\int\limits_0^1t^{b-1}{}_{r+1}F_{r}\left.\!\!\left(\begin{matrix}1-c+b,1-\f+b\\1-\f-\m+b\end{matrix}\right\vert
t\right)\frac{dt}{(1-t)^a}\\
=\frac{\Gamma(c)(\f-b)_\m}{\Gamma(b)\Gamma(c-b)(\f)_\m}\int\limits_0^1 \sum\limits_{n=0}^{\infty} \frac{(1-c+b)_n(1-\f+b)_n}{(1-\f-\m+b)_n n!} \frac{t^{b+n-1}}{(1-t)^a}dt
\\
=\frac{\Gamma(c)(\f-b)_\m}{\Gamma(b)\Gamma(c-b)(\f)_\m}\sum\limits_{n=0}^{\infty} \frac{(1-c+b)_n(1-\f+b)_n}{(1-\f-\m+b)_n n!}B(b+n;1-a)
\\
=\frac{\Gamma(c)\Gamma(1-a)(\f-b)_\m}{\Gamma(1+b-a)\Gamma(c-b)(\f)_\m}
{}_{r+2}F_{r+1}\left.\!\!\left(\begin{matrix}b,1-c+b,1-\f+b\\1+b-a,1-\f-\m+b\end{matrix}\right.\right),
\end{multline}
which is (\ref{eq:Kar4}) under additional conditions necessary for validity of (\ref{eq:int2}) and (\ref{eq:Ghiper}).
However,  both sides of (\ref{eq:Kar4}) are analytic in $a$, $b$, $c$, $f_i$, $i=1,\ldots r$ in the domain  $\Re(c-a-b-m)>0$, $c\notin-\N_0$ and $-f_i\notin\N_0$ for $i=1,\ldots,r$,
because
$$
(\f-b)_\m=(-1)^m\frac{\Gamma(1-\f+b)}{\Gamma(1-\f-\m+b)}
$$
so that the hypergeometric function on the right hand side is regularized.  Hence, (\ref{eq:Kar4}) remains true under the conditions stated in the theorem.  $\hfill\square$

\medskip

Note that another generalization of Karlsson's formula has been obtained by us in our recent work \cite[Corollaries~3,4]{KPdegMP} by replacing the unit argument in (\ref{eq:Karl1}) with an arbitrary argument $x$.

\section{Miller-Paris transformations and their consequences}

Recall that $\m=(m_1,\ldots,m_r)\in\N^r$, $m=m_1+m_2+\ldots+m_r$ and $\f=(f_1,\ldots,f_r)$.
Miller and Paris \cite[Theorem~3]{MP2013} established the following transformation:
\begin{equation}\label{eq:KRPTh1-1}
{}_{r+2}F_{r+1}\left.\!\!\!\left(\begin{matrix}a, b,
\f+\m\\c,\f\end{matrix}\right\vert x\right)
=(1-x)^{-a}{}_{m+2}F_{m+1}\left.\!\left(\begin{matrix}a,c-b-m, \zetta+1\\c,\zetta\end{matrix}\right\vert\frac{x}{x-1}\right),
\end{equation}
where $\zetta=\zetta(c,b,\f)=(\zeta_1,\ldots,\zeta_m)$ are the zeros of the characteristic polynomial
\begin{equation}\label{eq:Qm}
Q_m(t)=Q(b,c,\f,\m;t)=\frac{1}{(c-b-m)_{m}}\sum\limits_{k=0}^{m}(b)_kC_{k,r}(t)_{k}(c-b-m-t)_{m-k},
\end{equation}
with $C_{0,r}=1$, $C_{m,r}=1/(\f)_{\m}$ and
$$
C_{k,r}=\frac{1}{(\f)_{\m}}\sum\limits_{j=k}^m\sigma_j\mathbf{S}_j^{(k)}=\frac{(-1)^k}{k!}{}_{r+1}F_{r}\!\left(\begin{matrix}-k,\f+\m\\\f\end{matrix}\right).
$$
Here  $(\f)_{\m}=(f_1)_{m_1}\cdots(f_r)_{m_r}$ and  the numbers $\sigma_j$ ($0\leq{j}\leq{m}$) are defined by the generating function
\begin{equation}\label{eq:pol}
(f_1+x)_{m_1}\ldots(f_r+x)_{m_r}=\sum\limits_{j=0}^m\sigma_jx^j,
\end{equation}
where $\mathbf{S}_j^{(k)}$ is the Stirling number of the second kind.

Transformation (\ref{eq:KRPTh1-1}) fails when $c=b+p$, $p\in\N$, $1\leq p\leq{m}$. We explore in-depth what happens in this situation in our recent paper \cite{KPdegMP}.  Among other things we show that the characteristic polynomial (\ref{eq:Qm}) degenerates to
$$
R_{p-1}(a)=\sum\limits_{k=0}^m
{}_{r+2}F_{r+1}\left.\!\!\left(\begin{matrix}-k,\f+\m\\\f\end{matrix}\right.\right)\frac{(b)_k(a-k)_{p-1}}{k!}.
$$
In the theorem below we find a simpler expression for this polynomial which amounts to a generalization
of the summation formula  \cite[Corollary~3]{MP2012}.

\begin{theorem} \label{thm:determinate}
If $\Re(p+a-m-1)>0$, then
\begin{equation}
\sum\limits_{k=0}^m{}_{r+2}F_{r+1}\left.\!\!\left(\begin{matrix}-k,\f+\m\\\f\end{matrix}\right.\right)\frac{(b)_k(a-k)_{p-1}}{k!}=
\sum\limits_{q=0}^{p-1}\frac{(b)_q(\f-b-q)_\m}{(\f)_\m q!}(1-p)_q(b+q+a)_{p-1-q}.
\end{equation}
\end{theorem}
\textbf{Proof.} Making substitutions $a\to 1-a$ and $c\to{b+p}$ in (\ref{eq:sumMil}) and employing the elementary relation
$$
\frac{(1-a)_k}{(2-a-p)_k}=\frac{(a-k)_{p-1}}{(a)_{p-1}},
$$
we obtain
\begin{equation}\label{eq:sumMil6}
{}_{r+2}F_{r+1}\left.\!\!\left(\begin{matrix}1-a,b,\f+\m\\b+p,\f\end{matrix}\right.\right)=
\frac{\Gamma(b+p)\Gamma(p-1+a)}{\Gamma(b+p+a-1)\Gamma(p)(a)_{p-1}}\sum\limits_{k=0}^m
{}_{r+2}F_{r+1}\left.\!\!\left(\begin{matrix}-k,\f+\m\\\f\end{matrix}\right.\right)\frac{(b)_k(a-k)_{p-1}}{k!}.
\end{equation}
On the other hand, it follows from (\ref{eq:Karlson3}) that
\begin{equation}\label{eq:sumMil7}
{}_{r+2}F_{r+1}\left.\!\!\left(\begin{matrix}1-a,b,\f+\m\\b+p,\f\end{matrix}\right.\right)=
\frac{\Gamma(b+p)\Gamma(a)}{(\f)_\m}\sum\limits_{q=0}^{p-1}
\frac{(b)_q(\f-b-q)_\m}{(p-q-1)!(-q)_q\Gamma(a+q+b)}.
\end{equation}
Equating (\ref{eq:sumMil6}) and (\ref{eq:sumMil7}) and making simplifications in view of the identities
$$
\frac{(p-1)!}{(p-q-1)!(-q)_q}=\frac{(1-p)_q}{q!},
$$
$$
\frac{\Gamma(b+p+a-1)}{\Gamma(q+b+a)}=(b+q+a)_{p-q-1},
$$
we arrive at the statement of the theorem. $\hfill\square$

\medskip

In our next theorem we will write the summation formula (\ref{eq:sumMil}) in terms of the polynomial  $Q(b,c,\f,\m;t)$. This is not difficult to do, once the polynomial is transformed as follows:
\begin{equation}\label{eq:sumMil2}
Q(b,c,\f,\m;t)=\frac{(c-b-t-m)_m}{(c-b-m)_m}\sum\limits_{k=0}^m
{}_{r+1}F_{r}\left.\!\!\left(\begin{matrix}-k,\f+\m\\\f\end{matrix}\right.\right)
\frac{(t)_k(b)_k}{(1+t+b-c)_k k!}.
\end{equation}
However, we will give two different instructive proofs. The first shows how the Miller-Paris transformation formula (\ref{eq:KRPTh1-1}) leads directly to a summation formula. The second shows how this summation formula can be obtained from the integral representation (\ref{eq:int2}).

\begin{theorem} \label{thm:sumMil4}
Suppose $Q(b,c,\f,\m;t)$ is the polynomial defined in \emph{(\ref{eq:Qm})} and $\Re(c-a-b-m)>0$. Then
\begin{equation}\label{eq:sumMil4}
{}_{r+2}F_{r+1}\left.\!\!\left(\begin{matrix}a,b,\f+\m\\c,\f\end{matrix}\right.\right)=
\frac{\Gamma(c)\Gamma(c-a-b-m)}{\Gamma(c-a)\Gamma(c-b-m)}Q(b,c,\f,\m;a).
\end{equation}
\end{theorem}
\textbf{Proof~1.}  We will need the asymptotic formula \cite[(16.11.6)]{NIST}
\begin{equation}\label{eq:p+1Fp-asymp}
\frac{\Gamma(\a)}{\Gamma(\b)}{}_{p+1}F_{p}\left.\!\left(\begin{matrix}\a\\\b\end{matrix}\right\vert -z\right)=\sum_{k=1}^{p+1}\frac{\Gamma(a_k)\Gamma(\a_{[k]}-a_k)}{\Gamma(\b-a_k)}z^{-a_k}\bigl(1+\O(1/z)\bigr),~~\text{as}~|z|\to\infty,
\end{equation}
valid for  $|\arg(z)|<\pi$ if $a_j-a_k\notin\Z$ for all $k\ne{j}$.  Setting $-z=x/(x-1)$, we will have
\begin{multline}\label{eq:KRPTh1-11}
{}_{r+2}F_{r+1}\left.\!\!\left(\begin{matrix}a, b,\f+\m\\c,\f\end{matrix}\right\vert \frac{z}{1+z}\right)
=(1+z)^{a}{}_{m+2}F_{m+1}\left.\!\!\left(\begin{matrix}a,c-b-m, \zetta+1\\c,\zetta\end{matrix}\right\vert -z\right)
\\
=(1+z)^{a}\biggl\{z^{-a}\frac{\Gamma(c-b-m-a)\Gamma(\zetta+1-a)\Gamma(c)\Gamma(\zetta)}{\Gamma(c-b-m)\Gamma(\zetta-a)\Gamma(c-a)\Gamma(\zetta+1)}
\\
+C_0z^{-c+b+m}+C_1z^{-\zetta_1-1}+\ldots+C_m^{-\zetta_m-1}\biggr\}\bigl(1+\O(1/z)\bigr), \,\, z\to\infty.
\end{multline}
Here $C_0, C_1,\ldots,C_m$ are some constants. The roots of $Q(b,c,\f,\m;t)$ do not depend on $a$, so we can assume for the moment that  $\Re({-\zetta-1+a)}<0$. Further,  we need the condition $\Re(a-c+b+m)<0$ for the series defining the left hand side of (\ref{eq:sumMil4}) to converge. Letting $z\to\infty$ in (\ref{eq:KRPTh1-11}) we then get
\begin{equation}\label{eq:KRPTh1-14}
{}_{r+2}F_{r+1}\!\left(\begin{matrix}a, b,
\f+\m\\c,\f\end{matrix}\right)
=\frac{\Gamma(c-b-m-a)\Gamma(\zetta+1-a)\Gamma(c)\Gamma(\zetta)}{\Gamma(c-b-m)\Gamma(\zetta-a)\Gamma(c-a)\Gamma(\zetta+1)}.
\end{equation}
As $Q(b,c,\f,\m;0)=1$,
$$
\frac{\Gamma(\zetta+1-a)\Gamma(\zetta)}{\Gamma(\zetta-a)\Gamma(\zetta+1)}=\frac{(\zetta-a)_1}{(\zetta)_1}=Q(b,c,\f,\m;a)
$$
so that  (\ref{eq:KRPTh1-14}) takes the form (\ref{eq:sumMil4}). $\hfill\square$

\medskip

\textbf{Proof~2.} According to Theorem~\ref{thm:Karlson4} we have
\begin{equation*}
{}_{r+2}F_{r+1}\left.\!\!\left(\begin{matrix}a,b,\f+\m\\c,\f\end{matrix}\right.\right)=\frac{\Gamma(c)\Gamma(1-a)(\f-b)_\m}{\Gamma(1+b-a)\Gamma(c-b)(\f)_\m}
{}_{r+2}F_{r+1}\left.\!\!\left(\begin{matrix}b,1-c+b,1-\f+b\\1+b-a,1-\f-\m+b\end{matrix}\right.\right).
\end{equation*}
Representing the ${}_{r+2}F_{r+1}$ on the right hand side by the integral (\ref{eq:int2}) and using the form (\ref{eq:Ghiper}) for the density (the N{\o}rlund-Meijer $G$ function), similarly to the proof of  Theorem~\ref{thm:Karlson4}, in view of the identity (\ref{eq:Pocham}) we obtain:
\begin{equation}\label{eq:int5}
{}_{r+2}F_{r+1}\left.\!\!\left(\begin{matrix}a,b,\f+\m\\c,\f\end{matrix}\right.\right)
=\frac{\Gamma(c)}{\Gamma(c-b)\Gamma(b)}
\int\limits_{0}^{1}{}_{r+1}F_{r}\left.\!\!\left(\begin{matrix}a,\f+\m\\\f\end{matrix}\right\vert
t\right)\frac{t^{b-1}dt}{(1-t)^{1-c+b}}.
\end{equation}
Note that
$$
\frac{(\f+\m)_n}{(\f)_n}=\frac{1}{(\f)_m}\sum\limits_{k=0}^m\sigma_k n^k,
$$
where $\sigma_k$ is defined in (\ref{eq:pol}). Thus,
\begin{equation}\label{eq:int4}
\int\limits_{0}^{1}{}_{r+1}F_{r}\left.\!\!\left(\begin{matrix}a,\f+\m\\\f\end{matrix}\right\vert
t\right)\frac{t^{b-1}dt}{(1-t)^{1-c+b}}=\frac{1}{(\f)_m}\sum\limits_{k=0}^m\sigma_k\int\limits_0^1 \frac{t^{b-1}}{(1-t)^{1-c+b}}\sum\limits_{n=0}^{\infty}\frac{(a)_nt^n n^k}{n!}dt.
\end{equation}
The generating function of the Stirling numbers of the second kind is given by
$$
n^k=\sum\limits_{l=0}^k\mathbf{S}_k^{(l)}[n]_l,~~~\text{where}~~~[n]_l=n(n-1)\ldots(n-l+1),
$$
so that
\begin{multline*}
\sum\limits_{n=0}^{\infty} \frac{(a)_n t^n}{n!}n^k=\sum\limits_{n=0}^{\infty} \frac{(a)_n t^n}{n!}\sum\limits_{l=0}^k\mathbf{S}_k^{(l)}[n]_l=
\sum\limits_{l=0}^k\mathbf{S}_k^{(l)}\sum\limits_{n=l}^{\infty} \frac{(a)_n t^n}{(n-l)!}=\\
\sum\limits_{l=0}^k\mathbf{S}_k^{(l)}\sum\limits_{n=0}^{\infty} \frac{(a)_{n+l} t^{n+l}}{n!}=
\sum\limits_{l=0}^k\mathbf{S}_k^{(l)}\sum\limits_{n=0}^{\infty} \frac{(a)_l(a+l)_n t^{n+l}}{n!}=\sum\limits_{l=0}^k\mathbf{S}_k^{(l)}(a)_l\frac{t^l}{(1-t)^{a+l}}.
\end{multline*}
Substitution of this formula into (\ref{eq:int4}) yields:
\begin{equation}\label{eq:int6}
\int\limits_0^1 {}_{r+1}F_{r}\left.\!\!\left(\begin{matrix}a,\f+\m\\\f\end{matrix}\right\vert
t\right)\frac{t^{b-1}dt}{(1-t)^{1-c+b}}=\frac{1}{(\f)_m}\sum\limits_{k=0}^m\sigma_k\sum_{l=0}^{k}\mathbf{S}_k^{(l)}(a)_l\int\limits_0^1 \frac{t^{l+b-1}}{(1-t)^{a+l+1-c+b}}dt.
\end{equation}
Integrating, changing the order of summations and employing  (\ref{eq:int5}) we then obtain
$$
{}_{r+2}F_{r+1}\left.\!\!\left(\begin{matrix}a,b,\f+\m\\c,\f\end{matrix}\right.\right)=\frac{\Gamma(c)}{\Gamma(c-b)\Gamma(c-a)}
\sum_{l=0}^{m}\left(\frac{1}{(\f)_m}\sum\limits_{k=l}^m\sigma_k\mathbf{S}_k^{(l)}\right)(a)_l(b)_l\Gamma(c-a-b-l).
$$
In view of the definition of the polynomial $Q(b,c,\f,\m;t)$, simple algebra converts the right hand side of the above equality into the right hand side of (\ref{eq:sumMil4}). $\hfill\square$

\medskip

Note that the Miller-Paris transformation (\ref{eq:KRPTh1-1}) can be written as a transformation for the Meijer-N{\o}rlund $G$ function.  Indeed, suppose that $\f=(f_1,\ldots,f_r)$,  $b\neq f_j$, $1\leq j\leq r,$ $(c-b-m)_m\ne0$ and $\zetta^*$ is the vector of zeros of the polynomial $Q(b+1,c+1,\f+1,\m;t)$ defined in (\ref{eq:Qm}). According to (\ref{eq:KRPTh1-1}) we then have the identity
\begin{equation}\label{eq:FParis}
{}_{r+2}F_{r+1}\left.\!\!\left(\begin{matrix}a,b+1,\f+\m+1\\c+1,\f+1\end{matrix}\right\vert
x\right)=(1-x)^{-a}{}_{m+2}F_{m+1}\left.\!\!\left(\begin{matrix}a,c-b-m,\zetta^*+1\\c+1,\zetta^*\end{matrix}\right\vert
\frac{x}{x-1}\right)
\end{equation}
valid for any complex $a$.  At the same time, formula (\ref{eq:int2}) implies that
\begin{equation}\label{eq:F1}
{}_{r+2}F_{r+1}\left.\!\!\left(\begin{matrix}a,b+1,\f+\m+1\\c+1,\f+1\end{matrix}\right\vert
x\right)=\frac{\Gamma(c+1)}{\Gamma(b+1)(\f+1)_\m}\int\limits_0^{1}(1-xt)^{-a}G_{r+1,r+1}^{r+1,0}\!\left(\!t~\vline\begin{array}{l}c,\f\\
b,\f+\m\end{array}\!\!\right) dt.
\end{equation}

On the other hand, again using formula (\ref{eq:int2}) we obtain
\begin{multline}\label{eq:F2}
(1-x)^{-a}{}_{m+2}F_{m+1}\left.\!\!\left(\begin{matrix}a,c-b-m,\zetta^*+1\\c+1,\zetta^*\end{matrix}\right\vert
\frac{x}{x-1}\right)
\\
=(1-x)^{-a}\frac{\Gamma(c+1)}{\Gamma(c-b-m)(\zetta^*)_1}
\int\limits_0^{1}(1-(xt/(x-1))^{-a}
G_{m+1,m+1}^{m+1,0}\!\left(\!t~\vline\begin{array}{l}c,\zetta^*-1\\
c-b-m-1,\zetta^*\end{array}\!\!\right)dt
\\
=\frac{\Gamma(c+1)}{\Gamma(c-b-m)(\zetta^*)_1}\int\limits_0^{1}(1-xt)^{-a}G_{m+1,m+1}^{m+1,0}\!\left(\!1-t~\vline\begin{array}{l}c,\zetta^*-1\\
c-b-m-1,\zetta^*\end{array}\!\!\right)dt.
\end{multline}

Substitution of representations (\ref{eq:F1}) and (\ref{eq:F2}) into identity (\ref{eq:FParis}) yields a presumably new identity
\begin{equation}\label{eq:Gfunind}
G_{r+1,r+1}^{r+1,0}\!\left(\!t~\vline\begin{array}{l}c,\f\\
b,\f+\m\end{array}\!\!\right)=
\frac{\Gamma(b+1)(\f+1)_\m}{\Gamma(c-b-m)(\zetta^*)_1}
G_{m+1,m+1}^{m+1,0}\!\left(\!1-t~\vline\begin{array}{l}c,\zetta^*-1\\c-b-m-1,\zetta^*\end{array}\!\!\right).
\end{equation}
As $G$ function above reduces to a single hypergeometric function by (\ref{eq:Ghiper}), formula (\ref{eq:Gfunind})
can be reformulated in terms of hypergeometric functions. Indeed, applying (\ref{eq:Ghiper}) on both sides of (\ref{eq:Gfunind}), we obtain
\begin{equation*}
z^b{}_{r+1}F_{r}\left.\!\!\left(\begin{matrix}1-c+b,1-\f+b\\1-\f-\m+b\end{matrix}\right\vert
z\right)=K(1-z)^{c-b-m-1}{}_{m+1}F_{m}\left.\!\!\left(\begin{matrix}-b-m,1-\zetta^*+c-b-m\\-\zetta^*+c-b-m\end{matrix}\right\vert
1-z\right),
\end{equation*}
where
\begin{equation*}
K=\frac{(c-b-m)_m(\zetta^*-c+b+m)_1(\f+1)_\m}{(b+1)_m(\f-b)_\m(\zetta^*)_1}.
\end{equation*}
Combining (\ref{eq:Gfunind}) with the shifting property (\ref{eq:Gtimespower}) and the integral representation (\ref{eq:int2}), we arrive at the following statement.
\begin{theorem}\label{thm:product}
Suppose that $\Re(c-d-b-m)>0$, $\Re(a+b)>0$ and $(c-b-m)_m\ne0$.  Then the following identity holds:
\begin{multline}\label{eq:sum1}
{}_{r+2}F_{r+1}\left.\!\!\left(\begin{matrix}d,b+a,\f+\m+a\\c+a,\f+a\end{matrix}\right.\right)
\\
=\frac{\Gamma(c+a)\Gamma(b)(\f)_\m}{\Gamma(c)\Gamma(b+a)(\f+a)_\m}{}_{r+2}F_{r+1}\left.\!\!\left(\begin{matrix}d,b,\f+\m\\c,\f\end{matrix}\right.\right)
{}_{m+2}F_{m+1}\left.\!\!\left(\begin{matrix}-a,c-b-m-d,\zetta-d+1\\c-d,\zetta-d\end{matrix}\right.\right),
\end{multline}
where $\zetta$ is the vector formed by the zeros of the polynomial $Q(b,c,\f,\m;t)$ defined in \emph{(\ref{eq:Qm})}.
\end{theorem}
\textbf{Proof.} Applying the integral representation (\ref{eq:int2}), the shifting property (\ref{eq:Gtimespower}) and identity (\ref{eq:Gfunind}), we get:
\begin{multline*}
\frac{\Gamma(b+a)(\f+a)_\m}{\Gamma(c+a)}{}_{r+2}F_{r+1}\left.\!\!\left(\begin{matrix}d,b+a,\f+\m+a\\c+a,\f+a\end{matrix}\right.\right)
\\
=\int\limits_0^1(1-t)^{-d}G_{r+1,r+1}^{r+1,0}\!\left(\!t~\vline\begin{array}{l}c+a-1,\f+a-1\\b+a-1,\f+\m+a-1\end{array}\!\!\right)dt
\\
=\int\limits_0^1t^a(1-t)^{-d}G_{r+1,r+1}^{r+1,0}\!\left(\!t~\vline\begin{array}{l}c-1,\f-1\\b-1,\f+\m-1\end{array}\!\!\right)dt
\\
=\frac{\Gamma(b)(\f)_\m}{\Gamma(c-b-m)(\zetta)_1}\int\limits_0^1t^a(1-t)^{-d}G_{m+1,m+1}^{m+1,0}\!\left(\!1-t~\vline\begin{array}{l}c-1,\zetta-1\\
c-b-m-1,\zetta\end{array}\!\!\right)dt
\\
=\frac{\Gamma(b)(\f)_\m}{\Gamma(c-b-m)(\zetta)_1}\frac{\Gamma(c-b-m-d)(\zetta-d)_1}{\Gamma(c-d)}
{}_{m+2}F_{m+1}\left.\!\!\left(\begin{matrix}-a,c-b-m-d,\zetta-d+1\\c-d,\zetta-d\end{matrix}\right.\right).
\end{multline*}
Note that  $Q(b,c,\f,\m;0)=1$ and $(\zetta-d)_1/(\zetta)_1=Q(b,c,\f,\m;d)$. Then, using (\ref{eq:sumMil4}) to express $Q(b,c,\f,\m;d)$ we arrive at (\ref{eq:sum1}).
$\hfill\square$

In the final part of the paper we derive a recurrence relation connecting
$$
{}_{r+2}F_{r+1}\left.\!\!\left(\begin{matrix}a,b,\f+\m\\b+p,\f\end{matrix}\right.\right)
$$
with the functions of the type
$$
{}_{r+2}F_{r+1}\left.\!\!\left(\begin{matrix}a',b',\f'+\m\\b'+p-1,\f'\end{matrix}\right.\right)
$$
having a smaller negative parameter difference than the original function. Hence, iterative application of this recurrence gives yet another method to obtain an extension of Karlsson's theorem (\ref{eq:Karl1}).
\begin{theorem}\label{thm:recurs}
For integer $p\ge2$  and $\Re(p-a-m-2)>0$ the following recurrence relation is true:
\begin{multline}\label{eq:thm1}
{}_{r+2}F_{r+1}\left.\!\!\left(\begin{matrix}a,b,\f+\m\\b+p,\f\end{matrix}\right.\right)\!=\!\frac{(b+p-1)(p-a-1)}{(p-1)(b+p-a-1)}
{}_{r+2}F_{r+1}\left.\!\!\left(\begin{matrix}a,b,\f+\m\\b+p-1,\f\end{matrix}\right.\right)+\frac{ab}{(p-1)(b+p-a-1)}
\\
\times
\biggl\{{}_{r+2}F_{r+1}\left.\!\!\left(\begin{matrix}a+1,b+1,\f+\m\\b+p,\f\end{matrix}\right.\right)-
\frac{(\f+\m)_1}{(\f)_1}{}_{r+2}F_{r+1}\left.\!\!\left(\begin{matrix}a+1,b+1,\f+\m+1\\b+p,\f+1\end{matrix}\right.\right)\biggr\}.
\end{multline}
\end{theorem}
\textbf{Proof.} Denote
$$
S_p(a,b,\f,\m)=\sum\limits_{k=0}^m
{}_{r+1}F_{r}\left.\!\!\left(\begin{matrix}-k,\f+\m\\\f\end{matrix}\right.\right)
\frac{(a)_k(b)_k}{(a-p)_k k!}.
$$
Using the straightforward identities $(\alpha)_k/(\alpha-1)_k=(1+k/(\alpha-1))$ and $(\alpha)_{k+1}=\alpha(\alpha+1)_k$ we get
\begin{multline}\label{eq:rec1}
S_{p+1}(a,b,\f,\m)=\sum\limits_{k=0}^m
{}_{r+1}F_{r}\left.\!\!\left(\begin{matrix}-k,\f+\m\\\f\end{matrix}\right.\right)
\frac{(a)_k(b)_k(a-p)_k}{(a-p-1)_k (a-p)_k k!}
\\
=\sum\limits_{k=0}^m {}_{r+1}F_{r}\left.\!\!\left(\begin{matrix}-k,\f+\m\\\f\end{matrix}\right.\right)
\frac{(a)_k(b)_k}{(a-p)_k k!}\left(1+\frac{k}{a-p-1}\right)
\\
=S_{p}(a,b,\f,\m)+\frac{1}{a-p-1}\sum\limits_{k=1}^m{}_{r+1}F_{r}\left.\!\!\left(\begin{matrix}-k,\f+\m\\\f\end{matrix}\right.\right)\frac{(a)_k(b)_k}{(a-p)_k (k-1)!}
\\
=S_{p}(a,b,\f,\m)+\frac{ab}{(a-p)(a-p-1)}\sum\limits_{k=0}^{m-1}{}_{r+1}F_{r}\left.\!\!\left(\begin{matrix}-k-1,\f+\m\\\f\end{matrix}\right.\right)\frac{(a+1)_k(b+1)_k}{(a+1-p)_k(k)!}.
\end{multline}
Next, using the formula
\begin{equation}\label{eq:movie}
{}_{r+1}F_{r}\left.\!\!\left(\begin{matrix}-k-1,\a\\\b\end{matrix}\right.\right)
={}_{r+1}F_{r}\left.\!\!\left(\begin{matrix}-k,\a\\\b\end{matrix}\right.\right)-\frac{(\a)_1}{(\b)_1}
{}_{r+1}F_{r}\left.\!\!\left(\begin{matrix}-k,\a+1\\\b+1\end{matrix}\right.\right)
\end{equation}
and taking account of \cite[Corollary 1]{MP2012}
$$
{}_{r+1}F_{r}\left.\!\!\left(\begin{matrix}-m,\f+\m\\\f\end{matrix}\right.\right)=\frac{(-1)^m m!}{(\f)_\m},
$$
we conclude that
$$
{}_{r+1}F_{r}\left.\!\!\left(\begin{matrix}-m-1,\f+\m\\\f\end{matrix}\right.\right)=\frac{(-1)^m m!}{(\f)_\m}-\frac{(-1)^m m!(\f+\m)_1}{(\f)_1(\f+1)_\m}=0.
$$
 Hence, the summation on the right hand side of (\ref{eq:rec1}) can be extended to $k=m$. Application of (\ref{eq:movie}) on the right hand side of (\ref{eq:rec1}) transforms it as follows:
\begin{multline*}
S_{p+1}(a,b,\f,\m)=
S_{p}(a,b,\f,\m)+\frac{ab}{(a-p)(a-p-1)}\sum\limits_{k=0}^{m}
{}_{r+1}F_{r}\left.\!\!\left(\begin{matrix}-k,\f+\m\\\f\end{matrix}\right.\right)\frac{(a+1)_k(b+1)_k}{(a+1-p)_k (k)!}
\\
-\frac{ab(\f+\m)_1}{(a-p)(a-p-1)(\f)_1}\sum\limits_{k=0}^{m}
{}_{r+1}F_{r}\left.\!\!\left(\begin{matrix}-k,\f+\m+1\\\f+1\end{matrix}\right.\right)\frac{(a+1)_k(b+1)_k}{(a+1-p)_k (k)!},
\end{multline*}
or
\begin{multline}\label{eq:rec2}
S_{p+1}(a,b,\f,\m)=
S_{p}(a,b,\f,\m)+\frac{ab}{(a-p)(a-p-1)}S_{p}(a+1,b+1,\f,\m)
\\
-\frac{ab (\f+\m)_1}{(a-p)(a-p-1)(\f)_1}S_{p}(a+1,b+1,\f+1,\m).
\end{multline}
Further, according to  (\ref{eq:sumMil})
\begin{equation}\label{eq:relation}
{}_{r+2}F_{r+1}\left.\!\!\left(\begin{matrix}a,b,\f+\m\\b+p,\f\end{matrix}\right.\right)=\frac{B(b+p,p-a)}{B(p,b+p-a)}S_{p-1}(a,b,\f,\m).
\end{equation}
Writing  (\ref{eq:rec2}) for $S_{p-1}(a,b,\f,\m)$ we then obtain:
\begin{multline}\label{eq:rec3}
{}_{r+2}F_{r+1}\left.\!\!\left(\begin{matrix}a,b,\f+\m\\b+p,\f\end{matrix}\right.\right)=\frac{B(b+p,p-a)}{B(p,b+p-a)}
S_{p-2}(a,b,\f,\m)
\\
+\frac{B(b+p,p-a)}{B(p,b+p-a)}\left(\frac{abS_{p-2}(a+1,b+1,\f,\m)}{(a-p+2)(a-p+1)}-
\frac{ab (\f+\m)_1S_{p-2}(a+1,b+1,\f+1,\m)}{(a-p+2)(a-p+1)(\f)_1}\right).
\end{multline}
Applying (\ref{eq:relation}) once more we can transform the first term of the sum (\ref{eq:rec3}) as
\begin{multline*}
\frac{B(b+p,p-a)}{B(p,b+p-a)}
S_{p-2}(a,b,\f,\m)=\frac{B(b+p,p-a)}{B(p,b+p-a)}\frac{B(p-1,b+p-a-1)}{B(b+p-1,p-a-1)}
{}_{r+2}F_{r+1}\left.\!\!\left(\begin{matrix}a,b,\f+\m\\b+p-1,\f\end{matrix}\right.\right)
\\
=\frac{(b+p-1)(p-a-1)}{(p-1)(b+p-a-1)}
{}_{r+2}F_{r+1}\left.\!\!\left(\begin{matrix}a,b,\f+\m\\b+p-1,\f\end{matrix}\right.\right).
\end{multline*}
Transforming the other terms in (\ref{eq:rec3}) in a similar fashion yields (\ref{eq:thm1}). $\hfill\square$

\bigskip
\bigskip

\textbf{Acknowledgements.} We thank Michael Schlosser for pointing out to us the references \cite{Chu1,Chu2,Gasper,Haglund,Schlosser1,Schlosser2} and useful communications.
This research has been supported by the Russian Science Foundation under the project 14-11-00022.

\end{document}